\documentclass[11pt, a4paper]{article}
\usepackage{latexsym}
\usepackage{indentfirst}
\usepackage{graphicx}
\usepackage{amsmath}
\usepackage{amssymb}
\begin{document}
\title{Numerical investigation of the Bautin bifurcation in a
delay differential equation modeling leukemia}
\author{Anca Veronica Ion\\ \small{"Gh. Mihoc-C. Iacob" Institute of Mathematical Statistics
and Applied Mathematics}\\ \small{of the Romanian Academy, 13, Calea 13 Septembrie, 050711, }\\
\small{Bucharest, Romania; e-mail: anca-veronica.ion@ima.ro.} \\
Raluca Mihaela Georgescu \\\small{Faculty of Mathematics and Computer Sciences, University of Pite\c sti,}\\
\small{1, T\^ argu din Vale, 110040, Pite\c sti, Romania; e-mail: gemiral@yahoo.com.}}

\date{}
\maketitle

\begin{abstract}
In a previous work we investigated the existence of Hopf degenerate
bifurcation points for a differential delay equation modeling
leukemia and we actually found Hopf points of codimension two for
the considered problem. If around such a point we vary two
parameters (the considered problem has five parameters), then a
Bautin bifurcation should occur. In this work we chose a Hopf point
of codimension two for the considered problem and perform numerical
integration for parameters chosen in a neighborhood of the
bifurcation point parameters. The results show that, indeed, we have
a Bautin bifurcation in the chosen point.\\
\small{\textbf{Acknowledgement.} Work partially supported by Grant
11/05.06.2009 within the framework of the Russian Foundation for
Basic Research - Romanian Academy collaboration.}
\end{abstract}

\section{Introduction}
The considered equation was taken from  \cite{PM-M}, \cite{PM-B-M}
\begin{equation}\label{eq}
\dot{x}(t)=-\left[\frac{\beta_0}{1+x(t)^n}+\delta\right]x(t)+k\frac{\beta_0x(t-r)}{1+x(t-r)^n},
\end{equation}
and represents part of a model of periodic chronic myelogenous
leukemia.  The initial model consists of two delay equations, one
for the density of proliferating cells, $P$, and one for the density
of so-called ''resting'' cells, $N$. The latter equation is
independent, hence it can be studied independently (see \cite{PM-M},
\cite{PM-B-M}). Equation \eqref{eq} is this second equation, with
the unknown $N$ made dimensionless by dividing it by a quantity with
the same dimension. The parameters $\beta_0,\,n,\,\delta,\,k,\,r$
are positive real numbers. We do not insist here on their physical
significance since this is largely presented in \cite{PM-M},
\cite{PM-B-M}.
 Parameter $k$ is of the form $k=2e^{-\gamma r}$, with $\gamma$
 positive. We take here, as in \cite{AI-Hopf-deg}, $k$ as an independent
parameter, instead of $\gamma$. Note that, due to its definition,
$k<2$. We denote by $\alpha$ the vector of five parameters
($\alpha=(\beta_0,\,n,\,\delta,\,k,\,r)$).

The equilibrium points of the problem are, as can be easily seen
(\cite{PM-M}, \cite{PM-B-M})
$$x_1=0,\,\, x_2=x_2(\alpha)=(\frac{\beta_0}{\delta}(k-1)-1)^{1/n}.$$ The second
one is acceptable from the biological point of view if and only if
\begin{equation}\label{cond} \frac{\beta_0}{\delta}(k-1)-1
> 0,\end{equation}
condition that implies $k>1.$

The equilibrium point $x_2$ presents Hopf bifurcation for some
points in the parameter space \cite{PM-M}, \cite{PM-B-M},
\cite{AI-Hopf}. In \cite{AI-Hopf} we developed the apparatus for
investigating the normal form for a Hopf bifurcation point, by using
the center manifold theory (we needed a second order approximation
of the center manifold for this). In order to determine the normal
form, we computed the first Lyapunov coefficient, $l_1(\alpha)$.

In \cite{AI-Hopf-deg} we searched for points of degenerate Hopf
bifurcation, i.e. points $\alpha^*$ in the parameter space with
$l_1(\alpha^*)=0.$ The method used relied also on the center
manifold theory. We explored a quite extended zone of parameters,
having biological significance, and found that for $n=2$ and each
\break  $\beta_0\in \{0.5,1,1.5,2,2.5\},\, k\in\{1.1,1.2,...,1.9\}$,
values of $r$ and $\delta$ can be found, for which $l_1=0.$ We went
further, and for these points we computed the second Lyapunov
coefficient, by constructing a fourth order approximation of the
center manifold. For all the points $\alpha^*$  with
$l_1(\alpha^*)=0$ determined, we found  $l_2(\alpha^*)<0,$ and thus,
by the definition in \cite{Soto}, $(x_2,\,\alpha^*)$ represent Hopf
points of codimension two.

In the present paper we numerically investigate the occurrence of
Bautin bifurcation in one of the Hopf points of codimension two
found. In \cite{AI-Hopf-deg} we considered the restriction of the
problem to a two-dimensional center manifold, this restriction is a
two-dimensional problem and the theory of \cite{K} works for its
study. Some ideas concerning the restriction of the problem to the
center manifold are presented here, in Section 2.

In Section 3 we describe the Bautin bifurcation for two dimensional
dynamical systems, as it is presented in \cite{K}.

Then we chose a point in the parameter space for which $l_1=0$ and
we explore its neighborhood in order to see whether we regain the
bifurcation diagram of the Bautin bifurcation for our problem. For
this, we numerically integrate our problem for the chosen
parameters. As the graphs of the trajectories obtained by numerical
integration show, we have indeed a Bautin bifurcation in the chosen
Hopf point of codimension two (Section 4).

\section{The problem restricted to the center manifold}

In \cite{AI-Hopf-deg} we considered the nontrivial equilibrium point
$x_2$ and the linearized equation around this point, that is (see
also \cite{PM-M}, \cite{PM-B-M})
\begin{equation}\label{lineq}
\dot{z}(t)=-[B_1+\delta]z(t)+kB_1z(t-r),\end{equation} where
$z=x-x_2(\alpha),$ $B_1=\beta'(x_2)x_2+\beta(x_2)$, and
$\beta(x)=\displaystyle \frac{\beta_0}{1+x^n}.$

The characteristic equation corresponding to \eqref{lineq} is
\begin{equation}\label{chareq}
\lambda+\delta+B_1=kB_1e^{-\lambda r}.
\end{equation}

The eigenvalues depend on the vector of parameters, $\alpha$.

Assume that we have a point $\alpha^*$ in the parameters space such
that, for $\alpha$ in a neighborhood $U$ of $\alpha^*,$ there are
two eigenvalues, that we denote by $\lambda_{\alpha\,1,2}$ with the
property that all other eigenvalues have negative real part, and, at
$\alpha=\alpha^*,$ $\lambda_{\alpha^*,\,1,2}=\pm \omega^* i.$

The analysis in \cite{AI-Hopf} shows that
$\alpha^*=(\beta^*_0,\,n^*,\,\delta^*,\,k^*,\,r^*)$ satisfies the
above condition if and only if the relation
 \begin{equation}\label{surface-H}
r^*=\frac{\arccos((\delta^*+B^*_1)/(k^*B^*_1))}{\sqrt{(k^*B^*_1)^2-(\delta^*+B_1)^2}}\,,
\end{equation}
is satisfied, where $B_1^*$ is the value of $B_1$ at $\alpha^*.$

For $\alpha=\alpha^*$, a two-dimensional local invariant  manifold
(the local center manifold) exists and the reduction of the problem
to this manifold leads to the ordinary differential equation:

\begin{equation}\label{eq-series} \frac{du}{dt}=\omega^*i
u+\sum_{j+k\geq2}\frac{1}{j!k!}g_{jk}(\alpha^*)u^{j}\overline{u}^{k},
\end{equation}
where $u:\mathbb{R}\mapsto \mathbb{C}.$

The formalism for the construction of an approximation of the center
manifold and that for computing the coefficients $g_{jk}(\alpha^*)$
are fully presented in \cite{AI-Hopf-deg}. We  remind here only some
elements of that construction, that relies on the general ideas in
\cite{Far}, \cite{HaL}.

We considered the Banach space
$$\mathcal{B}=\left\{\psi:[-r,0]\mapsto \mathbb{R},\, \psi\,
\mathrm{is\, continuous\, on\,}[-r,0]  \right\},$$ and its
complexification, denoted $\mathcal{B}_{C}.$ We denoted by
$\mathcal{M}$ the subspace of $\mathcal{B}_{C}$ spanned by the two
eigenfunctions $\varphi_{\alpha^*\,1,2}(s)=e^{\pm\omega^*
is},\,s\in[-r,0],$ corresponding to the two eigenvalues
$\lambda_{\alpha^*\,1,2}$ and by $\mathcal{P}$ a projector defined
on $\mathcal{B}_{C},$ with values in $\mathcal{M}$. The local center
manifold is locally invariant, tangent to $\mathcal{M}$ in 0. It is
the graph of a smooth function
$w_{\alpha^*}:\mathcal{U}\subset\mathcal{M}\mapsto
(I-\mathcal{P})\mathcal{B}_{C}, \,$ ($\mathcal{U}$ is a neighborhood
of 0 in $\mathcal{M}$)  that satisfies $w_{\alpha^*}(0)=0.$ Thus a
``point'' $\phi$ on the center manifold has the form $\phi
=z\varphi_{\alpha^*\,1}+\overline{z}\varphi_{\alpha^*\,2}+w_{\alpha^*}(z\varphi_{\alpha^*\,1}+\overline{z}\varphi_{\alpha^*\,2}).$

For an initial condition $\phi$ on the center manifold, the solution
$x(\cdot)$ of equation \eqref{eq} satisfies
\[x_t=u(t)\varphi_{\alpha^*\,1}+\overline{u}(t)\varphi_{\alpha^*\,2}+w_{\alpha^*}(u(t)\varphi_{\alpha^*\,1}+\overline{u}(t)\varphi_{\alpha^*\,2}),
\]
where $x_t\in \mathcal{B},$ is defined by
$x_t(s)=x(t+s),\,s\in[-r,0],$  $u(\cdot)$ is the solution of
equation \eqref{eq-series} with the initial condition $u(0)=u_0,$
and $\mathcal{P}\phi=u_0\varphi_1+\overline{u}_0\varphi_2$.

When $\alpha\in U$, the two eigenvalues
$\lambda_{\alpha,\,1,2}=\mu(\alpha)\pm i \omega(\alpha)$ may have
positive or negative real part.
 In each of these situations,  there still is a
two-dimensional local invariant manifold, a local unstable manifold
when $Re \lambda_{\alpha\,1,2}>0$, and a submanifold of the local
stable manifold for $Re \lambda_{\alpha\,1,2}<0$. This latter case
can be argued with the ideas of \cite{foliations}, adapted to the
more simple case considered by us. Hence the solution of our problem
has in this case also a representation of the form
\[x_t=u(t)\varphi_{\alpha\,1}+\overline{u}(t)\varphi_{\alpha\,2}+w_{\alpha}(u(t)\varphi_{\alpha\,1}+\overline{u}(t)\varphi_{\alpha\,2}),
\]
where u satisfies an equation of the form
\begin{equation}\label{eq-series1} \frac{du}{dt}=(\mu+\omega i)
u+\sum_{j+k\geq2}\frac{1}{j!k!}g_{\alpha\,jk}(\alpha)u^{j}\overline{u}^{k},
\end{equation}
and $w_{\alpha}$ is the function whose graph is the local invariant
manifold.

 The above considerations show that the  problem \eqref{eq-series} (respectively \eqref{eq-series1}) presents,
for some initial value, a periodic solution iff the corresponding
solution $x(t)$ of \eqref{eq} is periodic. Also, a solution of
\eqref{eq-series} (respectively \eqref{eq-series1}) spirals towards
0 (or from 0), iff the corresponding solution of \eqref{eq} spirals
towards (respectively from 0). Hence the study of equations
\eqref{eq-series}, \eqref{eq-series1} from the point of view of the
Hopf or Bautin bifurcation leads to complete conclusions concerning
these bifurcations for the problem \eqref{eq}.

\section{Bautin bifurcation for planar systems \cite{K}}
Consider a system of two ODEs, that can be written as a single
complex equation as
\begin{equation}\label{K1}\dot{z}=\lambda(\alpha)z+\sum_{j+k\geq2}\frac{1}{j!k!}g_{jk}(\alpha)z^{j}\overline{z}^{k},
\end{equation}
where $\alpha=(\alpha_1,\alpha_2)\in \mathbb{R}^2.$

In the hypotheses that a certain value $\alpha_0$ of $\alpha$ exists
such that:
\begin{itemize}
    \item $\lambda(\alpha_0)= i \omega_0$,
    \item $l_1(\alpha_0)= 0 $,
    \item $l_2(\alpha_0)\neq 0 $,
    \item the map  $\alpha \rightarrow
(\mu_1(\alpha),\,\mu_2(\alpha))$, where
$\mu_1(\alpha)=\frac{\mu(\alpha)}{\omega(\alpha)},\,\mu_2(\alpha)=l_1(\alpha)$
is regular at $\alpha_0$,
\end{itemize}
equation \eqref{K1} may brought by several transform of functions
and of parameters to the form:
\begin{equation}\label{K2}
\dot{u}=(\mu_1(\alpha)+i)u+\mu_2(\alpha)u|u|^{2}+L_2(\alpha)u|u|^4+O(|u|^5),
\end{equation}
where $u:\mathbb{R}\mapsto \mathbb{C}$ and
$L_2(\mu)=l_2(\alpha(\mu))$.

Moreover, in \cite{K} is proved that eq. \eqref{K2} is locally
topologically equivalent with
\begin{equation}\label{K3}
\dot{u}=(b_1+i)u+b_2 u|u|^{2}+su|u|^4,
\end{equation}
where $b_1=\mu_1$,  $b_2=\sqrt{|L_2(\mu)|}\mu_2$ and $s$ is the
signature of $l_2(\alpha_0).$

In order to describe the phase portrait for the parameters varying
around the point $\alpha_0$ (equivalent to $(b_1,b_2)=(0,0)$) it is
useful to consider the polar form of the above equation:\\
 $$\dot{\rho}=\rho(b_1+b_2\rho^2+s\rho^4),$$
 $$\dot{\theta}=1.$$\\
The limit cycles are obtained by solving the equation:
$$b_1+b_2\rho^2+\rho^4=0,
$$
and, by studying the number of its solutions as function of $b_1,
b_2$, the bifurcation diagram of the Bautin bifurcation is obtained.

We reproduce in Fig. 1 the bifurcation diagram for the case $s=-1$
since for all our Bautin bifurcation points found the second
Lyapunov coefficient is negative.

\begin{figure}\centering
\includegraphics[width=0.50\linewidth]{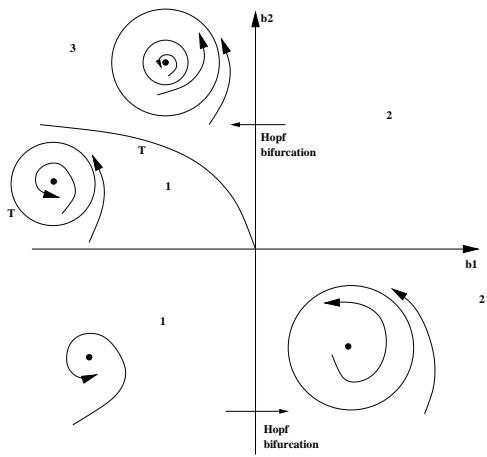}\caption{\small{Bautin bifurcation diagram for $\mathrm{sign} l_2(\alpha_0)<0$.}}
\end{figure}
We see that in a neighborhood of the origin in the plane of the
parameters $(b_1, \,b_2)$ the phase portrait in a neighborhood of
$u=0$ has very different aspects. These are described in \cite{K},
but for the sake of completeness, we point out a few ideas here.  In
the zone 1 of the Bautin bifurcation diagram, that lies between the
$b_2<0$ part of the axis $b_1=0$ and the curve $T$,
 the point $u=0$ is an attractive focus; when we cross the axis
 $b_1=0$ entering in the zone 2 (the half-plane $b_1>0$), a supercritical Hopf bifurcation takes place
 and a stable limit cycle occurs, while the point $u=0$ loses
 stability; then, starting from the first quadrant, when we cross
 the  axis $b_1=0$ to arrive in the zone 3 (lying between the $b_2>0$ part of the axis $b_1=0$ and the
 curve $T$), a subcritical Hopf
 bifurcation takes place, and an unstable limit cycle is born, in
 the interior of that previously formed. Then, for the parameters $(b_1,\,b_2)$
 on the curve $T$, the two cycle collide in a single cycle, that is
 repulsive on its interior side and attractive on its exterior side, and after crossing the curve
 $T$, arriving in the zone 1 again,
 the two cycles disappear.

Hence, the most interesting feature of this bifurcation is the
presence of two limit cycles one inside the other for the parameters
lying between the axis $b_1=0$ and the curve T. The exterior cycle
is stable (attractive) while the interior one is unstable
(repulsive).

\subsection{Numerical confirmation of Bautin bifurcation}

We have chosen the case $n^*=2,\,\beta^*_0=2.5,\,k^*=1.01,$ that is
not among those found in \cite{AI-Hopf-deg}. We took this case
hoping to have a small $r$ at the Bautin bifurcation and intending
to have the second Lyapunov coefficient not very close to zero (this
choice is justified by the table contained in Fig. 6 of
\cite{AI-Hopf-deg}).

By using the methods  presented in  \cite{AI-Hopf-deg}, we find that at\\
$r^*_0=5.301432998,\,\delta^*_0=0.0023073665$, we have $l_1=0,$
while $l_2=-0.0662.$

In order to see if we can regain the bifurcation diagram above for
our problem, we performed numerical integrations of the delay
differential equation \eqref{eq} for the above values of
$n,\,\beta_0,\,k$ and values of $r,\,\delta$ around
$r^*,\,\delta^*_0.$ We used the routine dde23 of Matlab.

In Fig. 2 we see, in the plane $(r,\delta)$, the curve of points
where $\omega=0$, the point $B$ where $l_1=0$, and the points where
we performed the numerical integration.

\begin{figure}\centering
\includegraphics[width=0.44\linewidth]{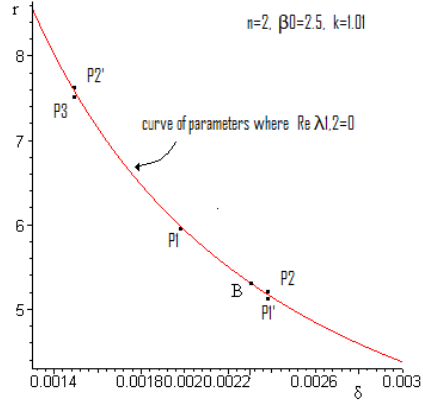}
\caption{\small{Curve of points $(\delta,r)$ where $Re
\lambda_{1,2}=0$, for $n=2,\,\beta_0=2.5,\,k=1.01$, the point of
Bautin bifurcation (point $B$), and the points chosen for numerical
integration.}}
\end{figure}

In order to  put into light the behavior of the solution, that is
qualitatively described in the bifurcation diagram, for a chosen
point in the parameters space, we have to take several initial
functions situated at different distances from the equilibrium
point. We took the initial function $\phi$ of the form
$\phi=x_2+c\,e^{\mu s}\cos(\omega s),$ where $\lambda_{\alpha
1,2}=\mu\pm i\omega$ are the eigenvalues of the linearized problem
at the chosen parameters, and $c$ is a new parameter, that we vary.

The results of the integrations, for each of the considered points
and for some choices of $c$ are represented vs time, but also in
$x(t),\,\dot{x}(t)$ plots.

We consider first the solutions for two points $P_1,\,P_1'$ in the
zone of the $(\delta,\,r)$ plane, corresponding to the zone 1 of the
bifurcation diagram. More precisely the point $P_1$ has the
coordinates $\delta=0.002,\;r=5.93,$ while for $P_1',\,
\delta=0.0024,\,r=5.14$ (see Fig. 2).
\begin{figure}\centering
\includegraphics[width=0.49\linewidth]{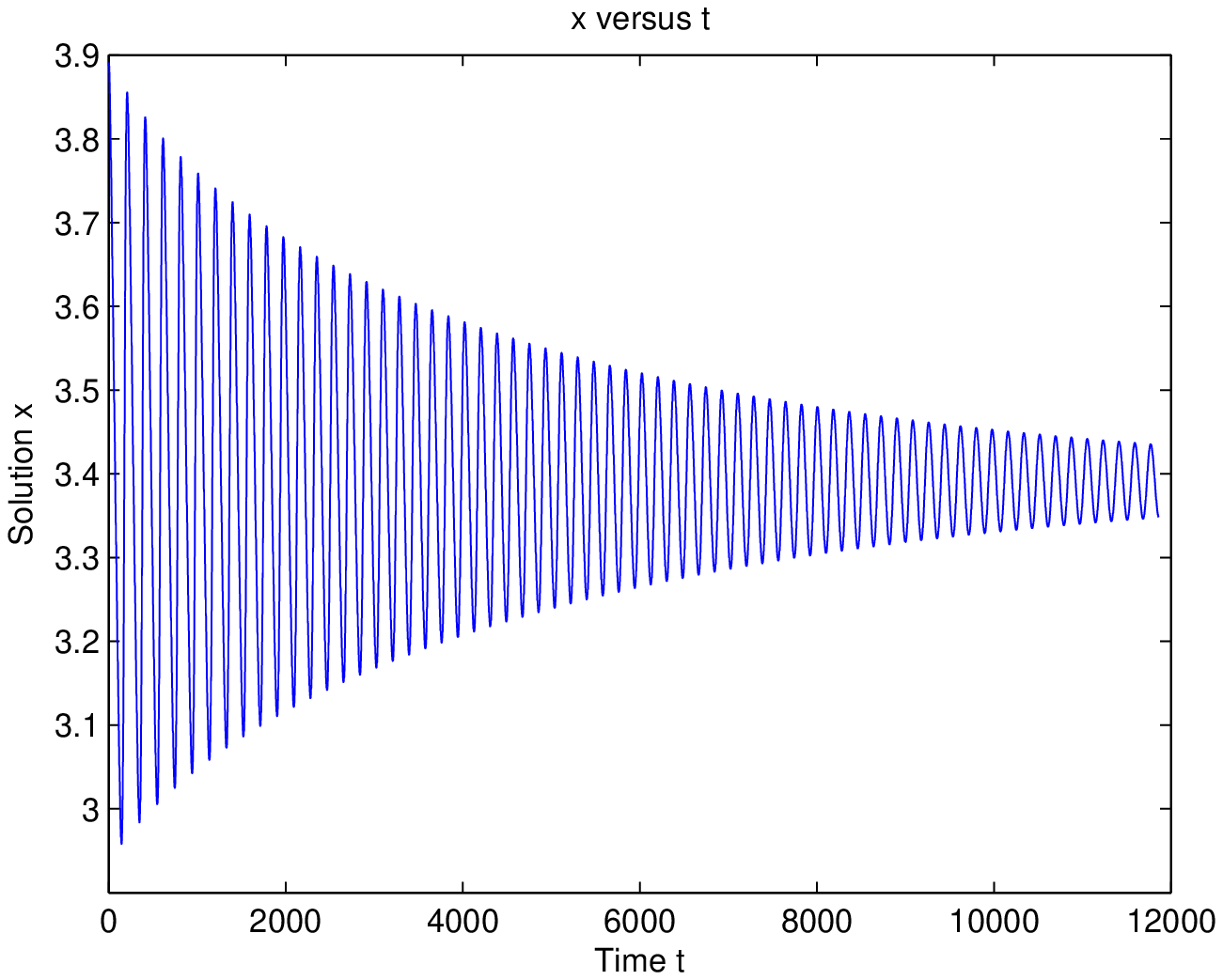}\includegraphics[width=0.49\linewidth]{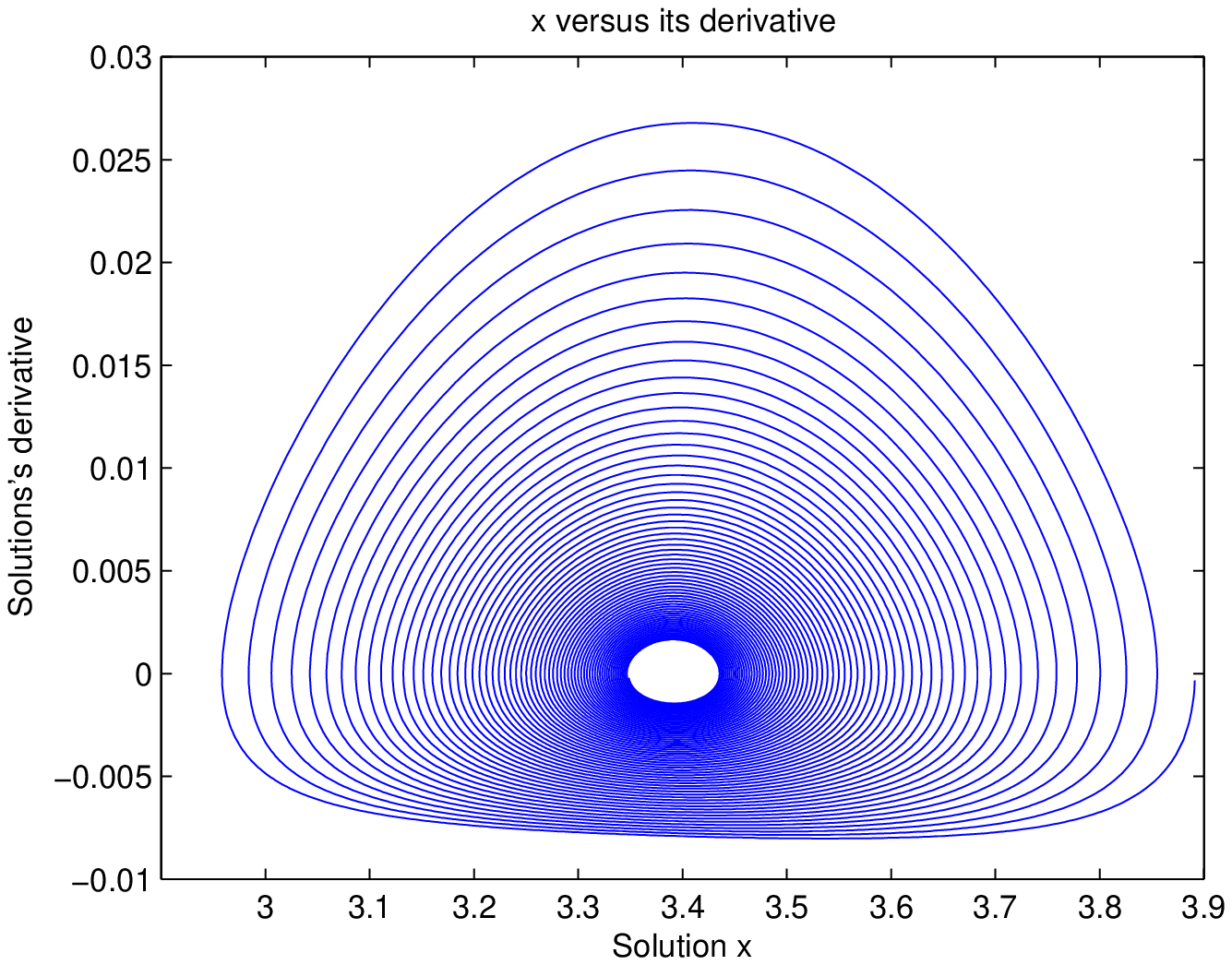}
\caption{\small{Phase portrait for the parameters in the point
$P_1$, for $c=0.5$.}}
\end{figure}
\begin{figure}\centering
\includegraphics[width=0.49\linewidth]{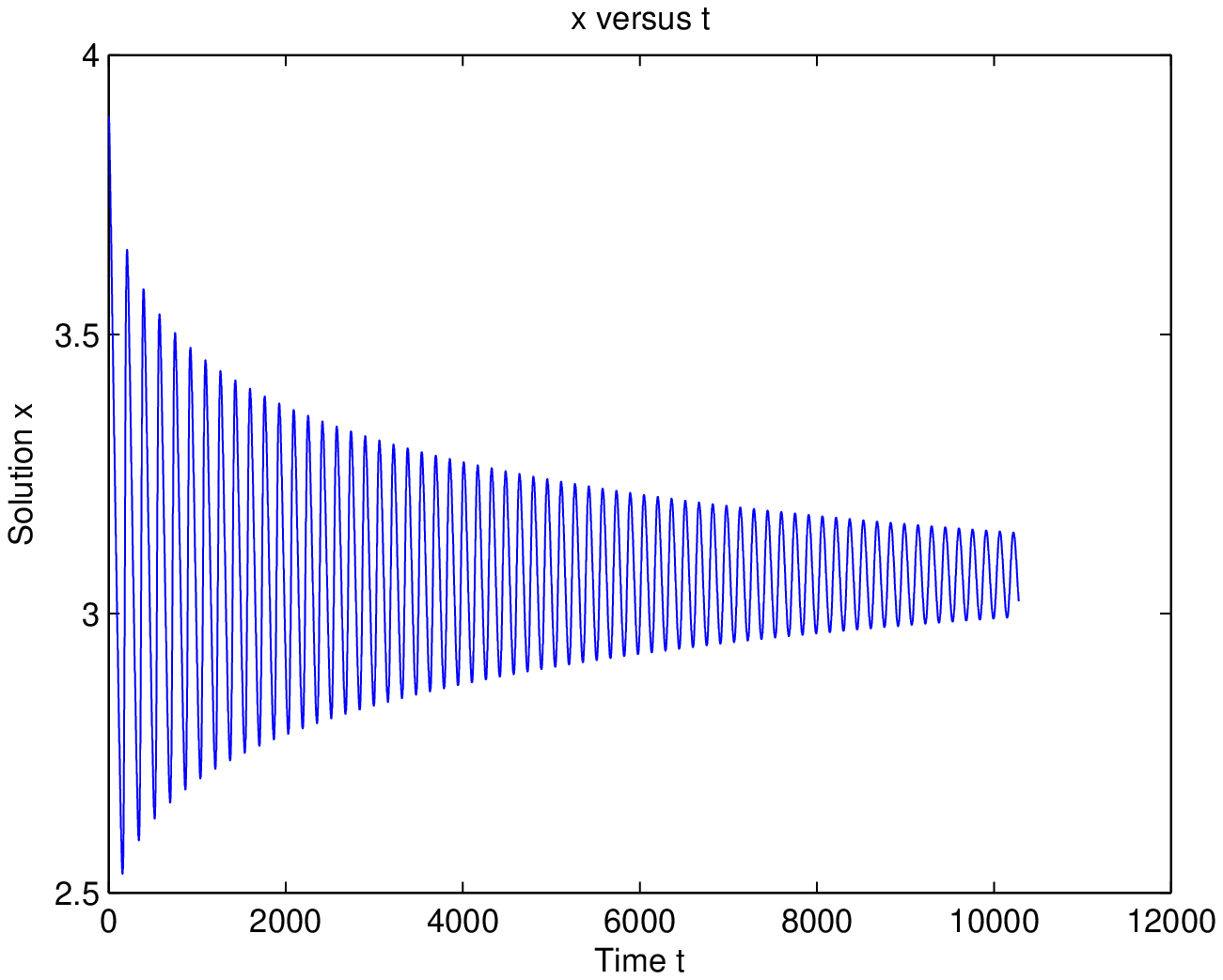}\includegraphics[width=0.49\linewidth]{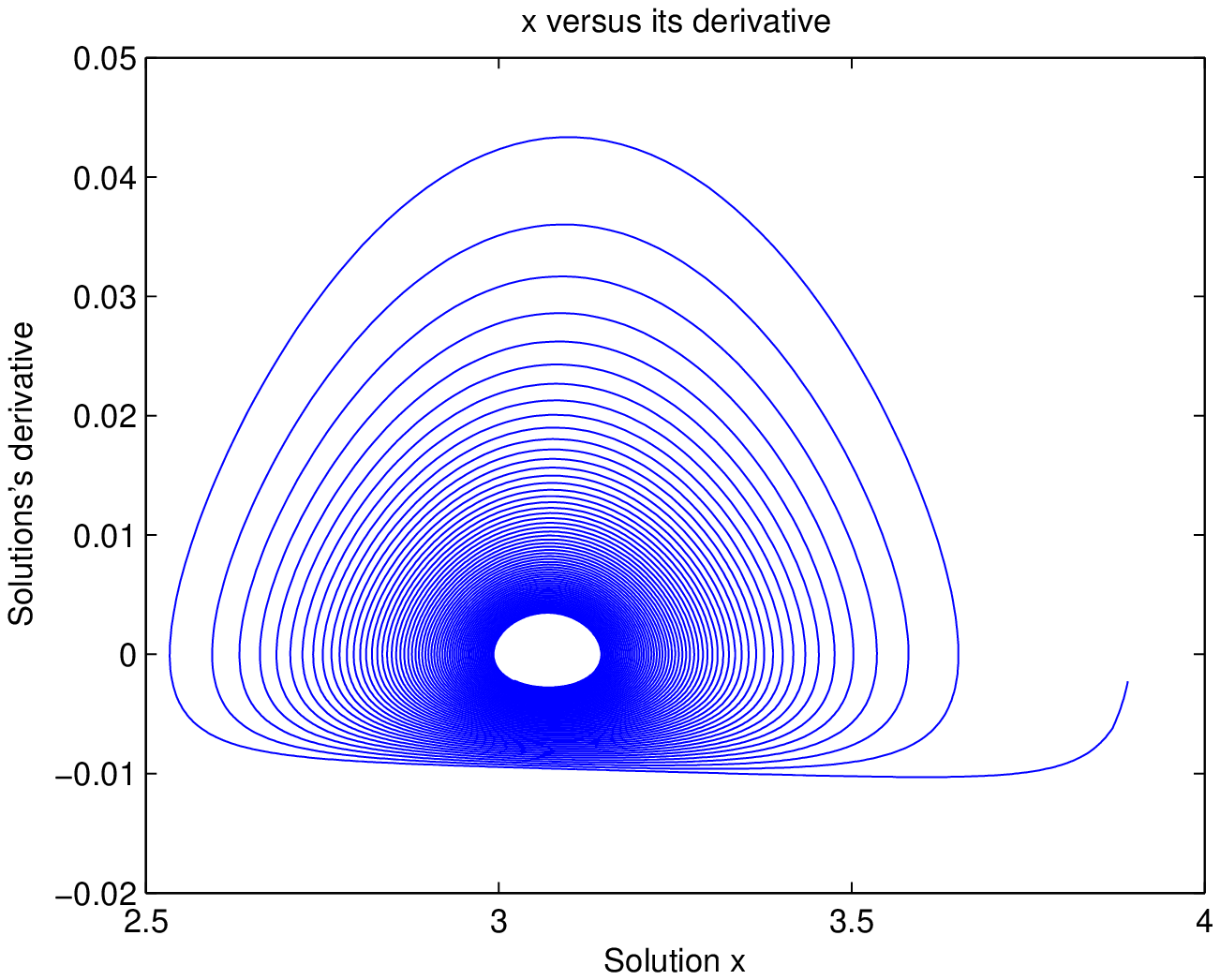}
\caption{\small{Phase portrait for the parameters in the point
$P_1'$, for $c=0.5$.}}
\end{figure}

The behavior of the solution for these two points is that
corresponding to an attracting focus. Since we obtained
qualitatively the same image for several values of $c$, we represent
only one of these, that for $c=0.5$, in Figs. 3 and 4.

The two points considered next, $P_2\,(\delta=0.0024,\,r=5.2)$ and
$P_2'\,(\delta=0.0015,\,r=7.56)$ are situated in the zone
corresponding to zone 2 of the bifurcation diagram (see Fig. 2).

\begin{figure}\centering
\includegraphics[width=0.44\linewidth]{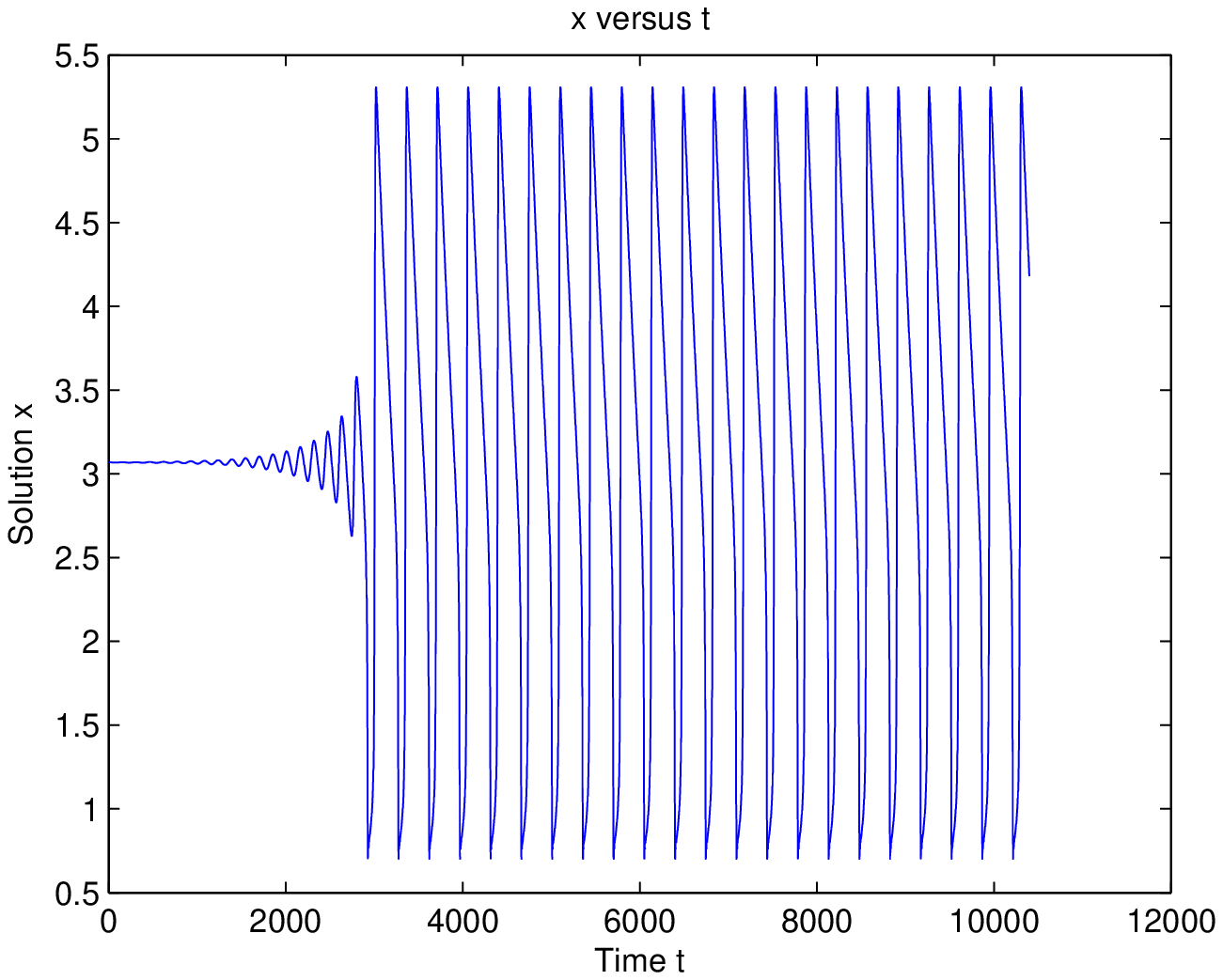}\includegraphics[width=0.49\linewidth]{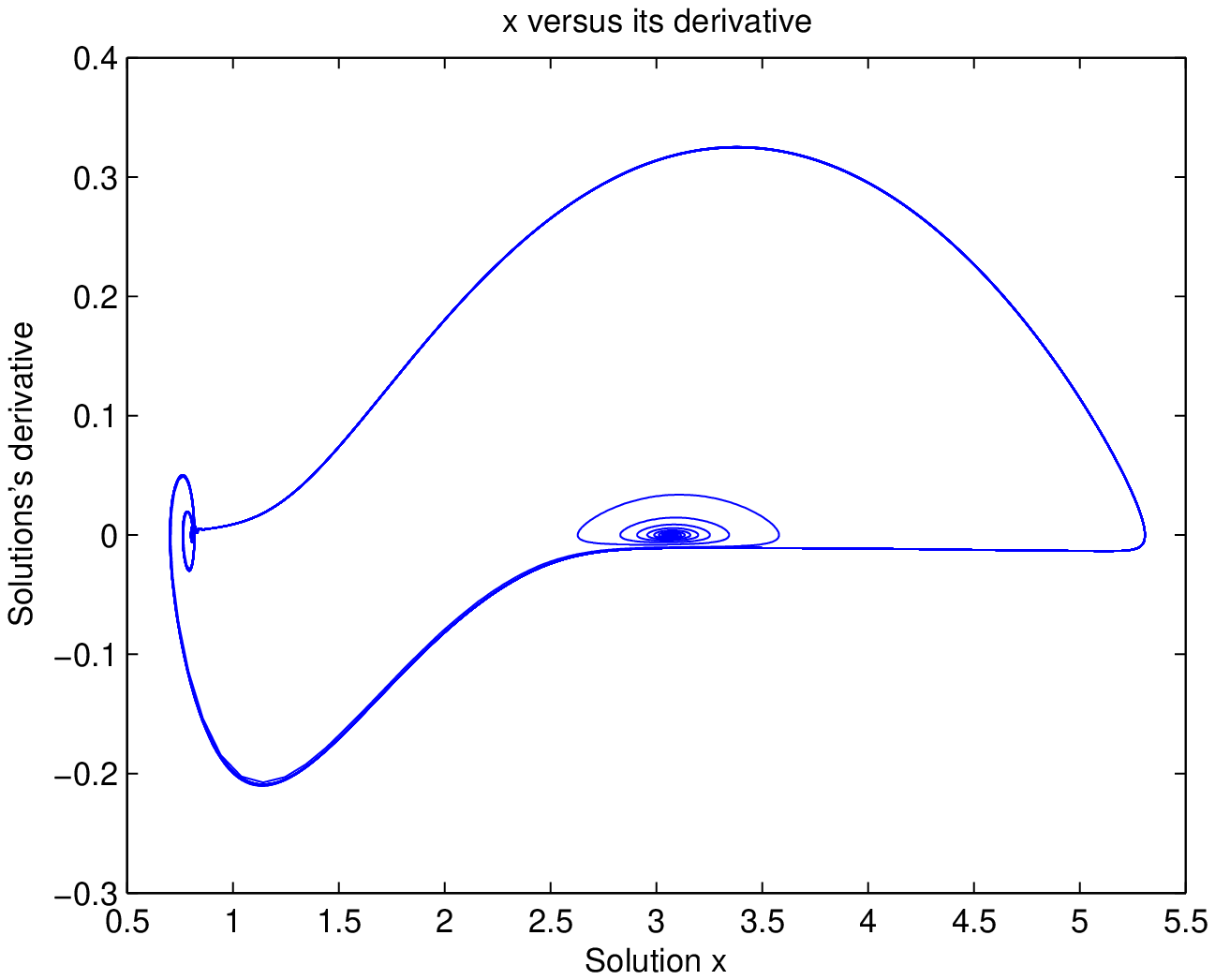}
\caption{\small{Phase portrait for the parameters in the point
$P_2$, for $c=0.001$.}}
\end{figure}
\begin{figure}\centering
\includegraphics[width=0.44\linewidth]{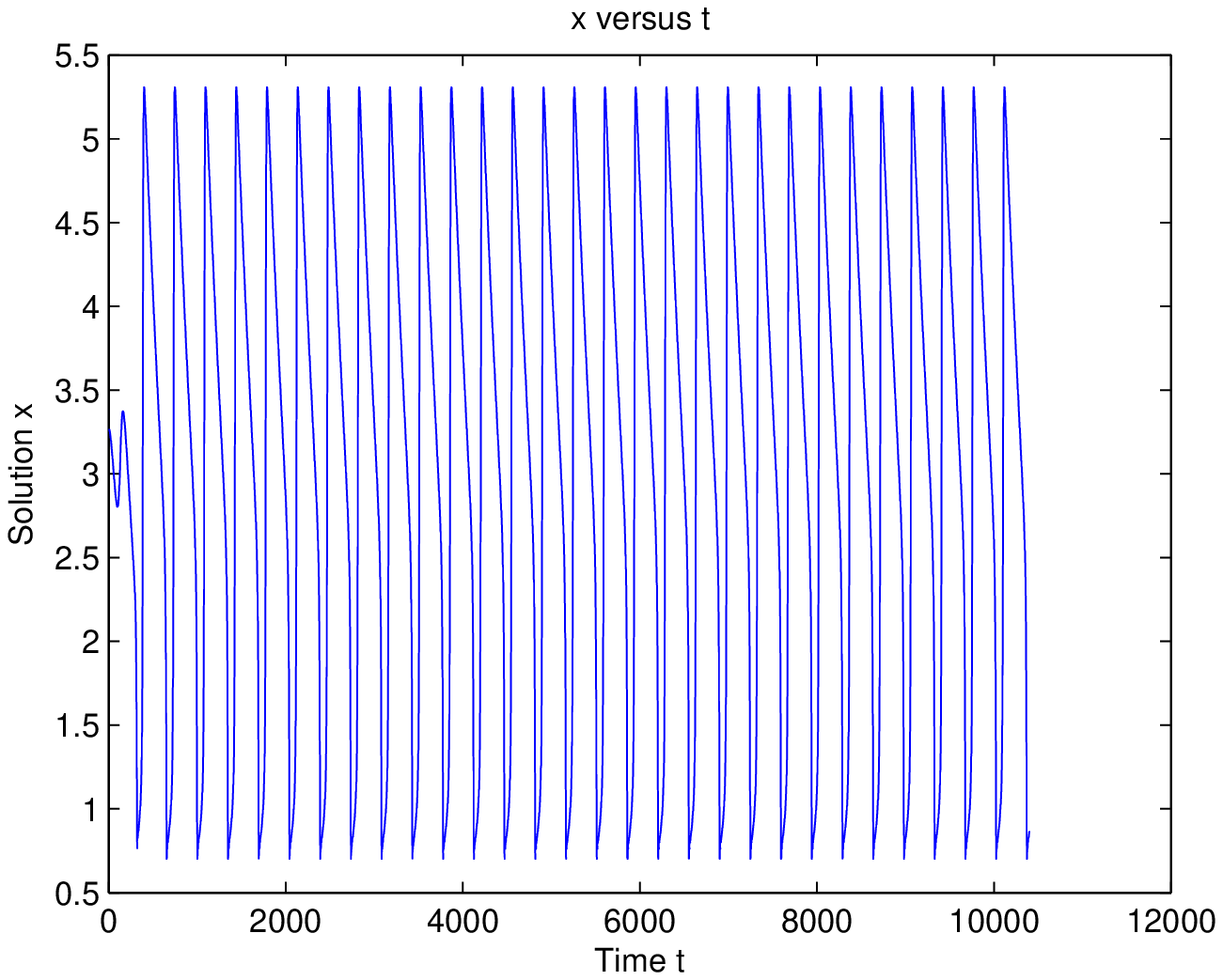}\includegraphics[width=0.49\linewidth]{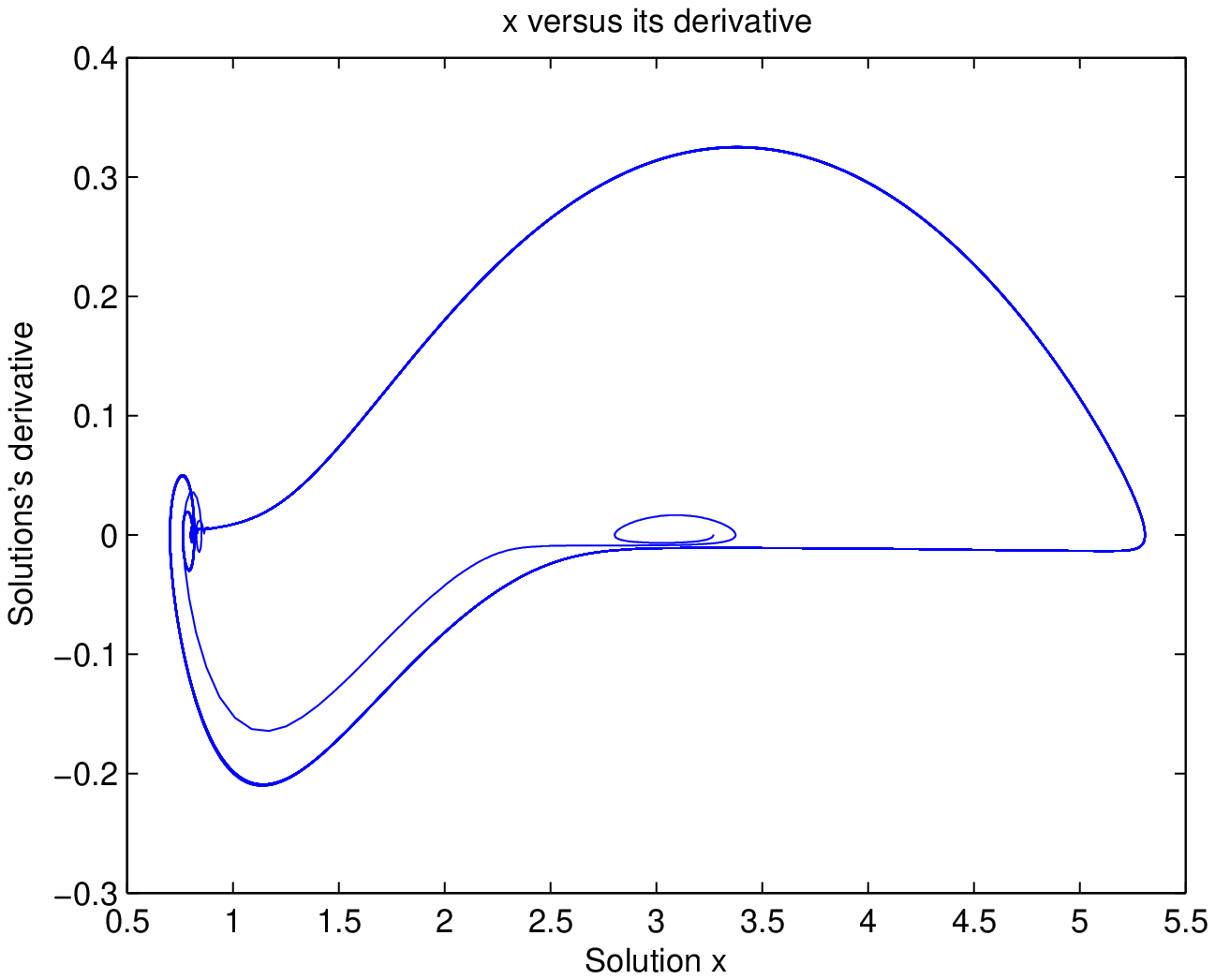}
\caption{\small{Phase portrait for the parameters in the point
$P_2$, for $c=0.2$.}}
\end{figure}
We see in Figs. 5 - 6 that a stable limit cycle occurs by Hopf
bifurcation. In these two figures we present the behavior of the
solution corresponding to $P_2$ and $c=0.001,$ respectively $c=0.2$.
In Figs. 7 and 8 we present the behavior of the solution
corresponding to $P_2'$ and $c=0.1,$ respectively $c=5$.
\begin{figure}\centering
\includegraphics[width=0.44\linewidth]{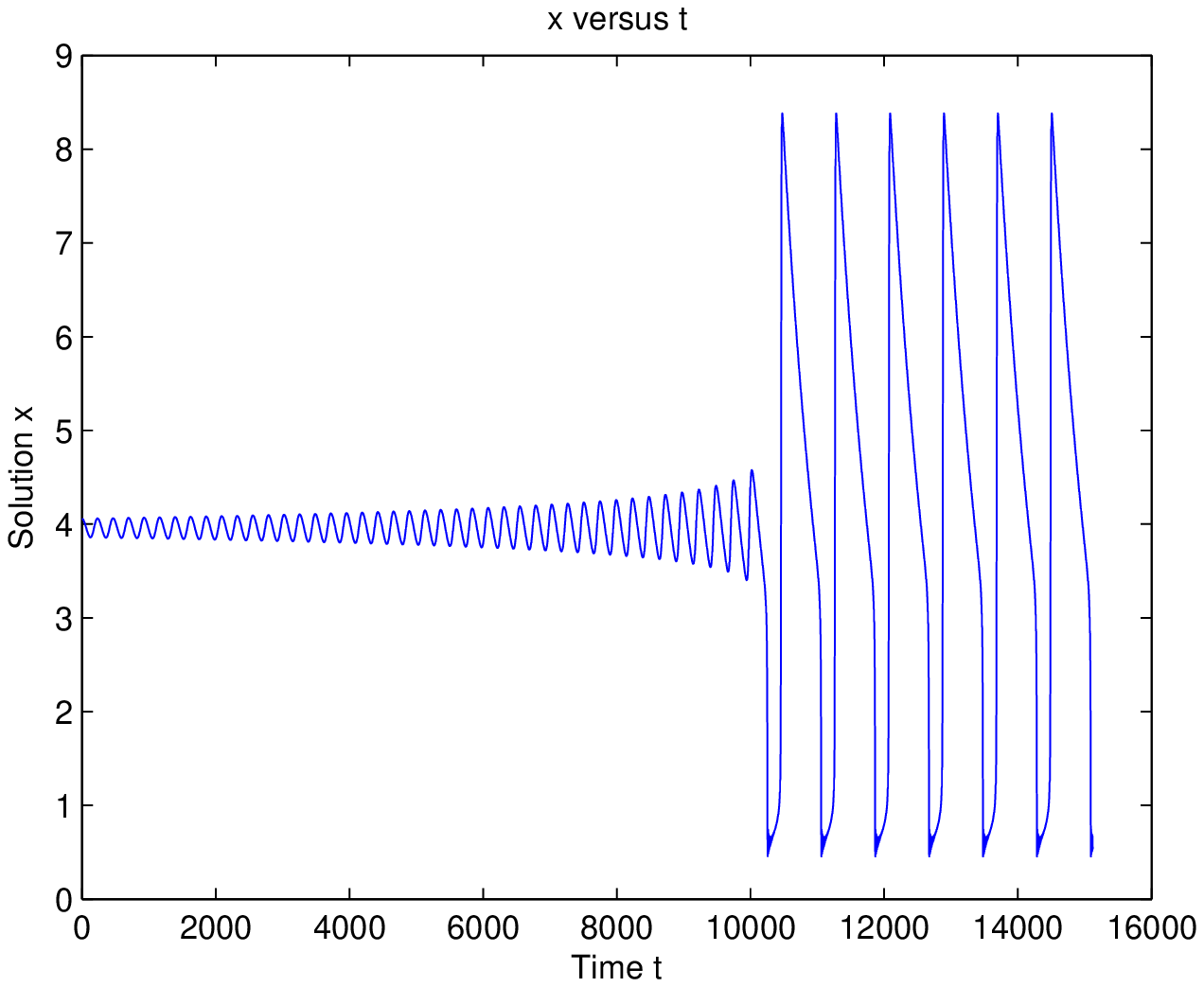}\includegraphics[width=0.49\linewidth]{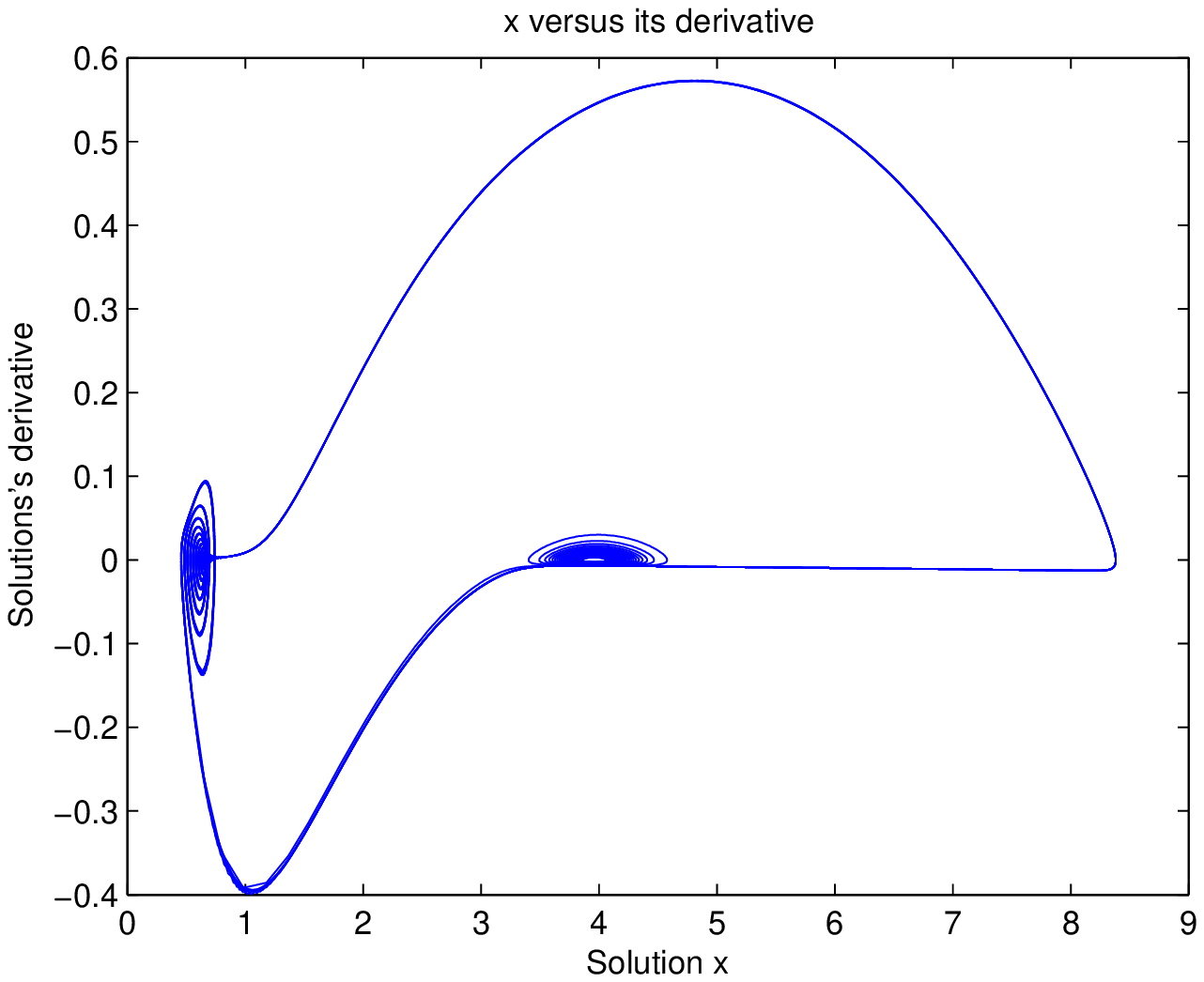}
\caption{\small{Phase portrait for the parameters in the point
$P_2'$, for $c=0.1$.}}
\end{figure}
\begin{figure}\centering
\includegraphics[width=0.44\linewidth]{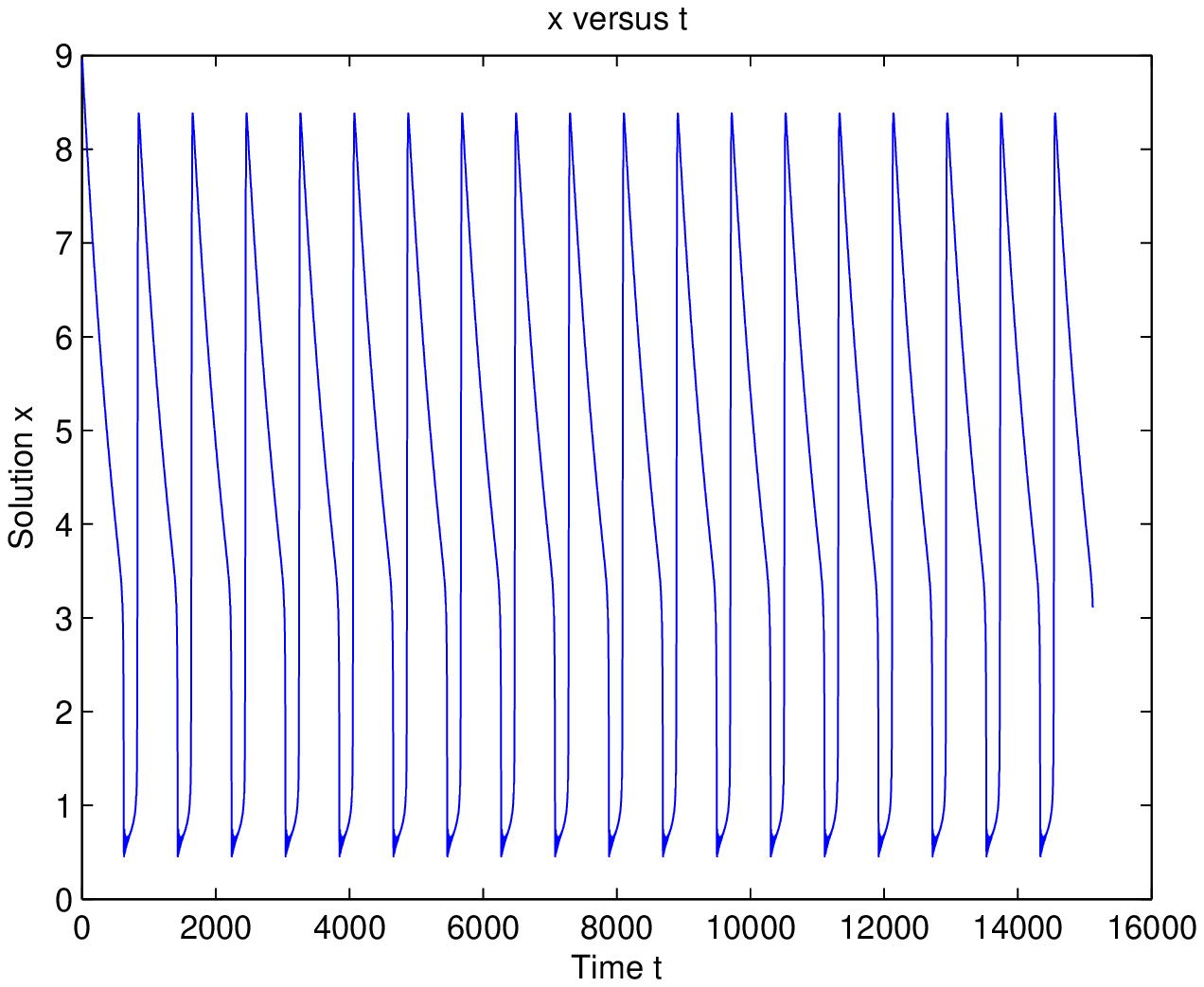}\includegraphics[width=0.49\linewidth]{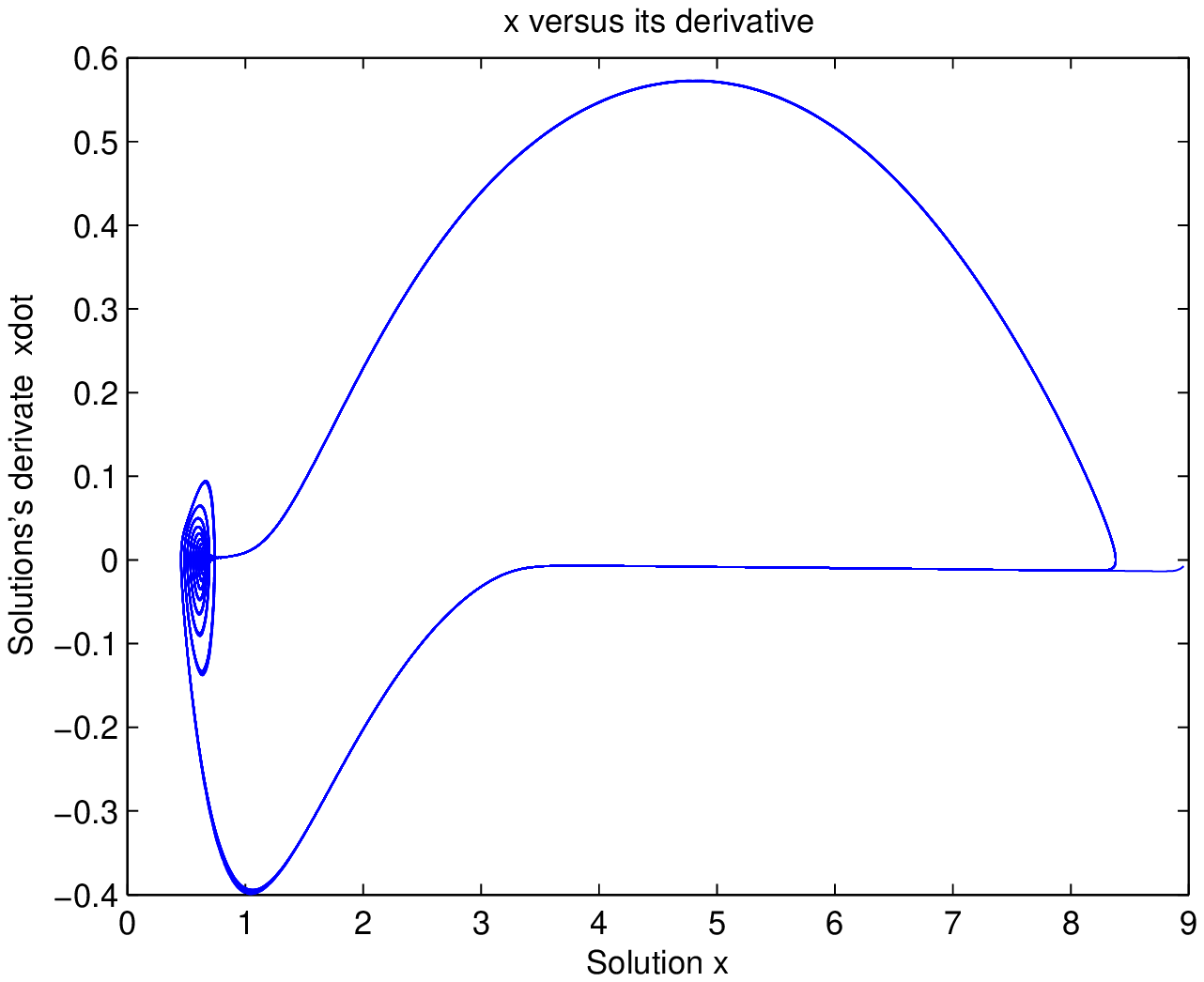}
\caption{\small{Phase portrait for the parameters in the point
$P_2'$, for $c=5$.}}
\end{figure}
We remark that at each oscillation, on this limit cycle, at the end
of a descending branch, in the $x$ versus $t$ representations, a
small superposed oscillation occurs, that produces a little spiral
in the left of the $x$ versus $\dot{x}$ representation. The moment
when the solutions enters on this limit cycle depends on the
distance between the initial function and $x_2$.
\begin{figure}\centering
\includegraphics[width=0.44\linewidth]{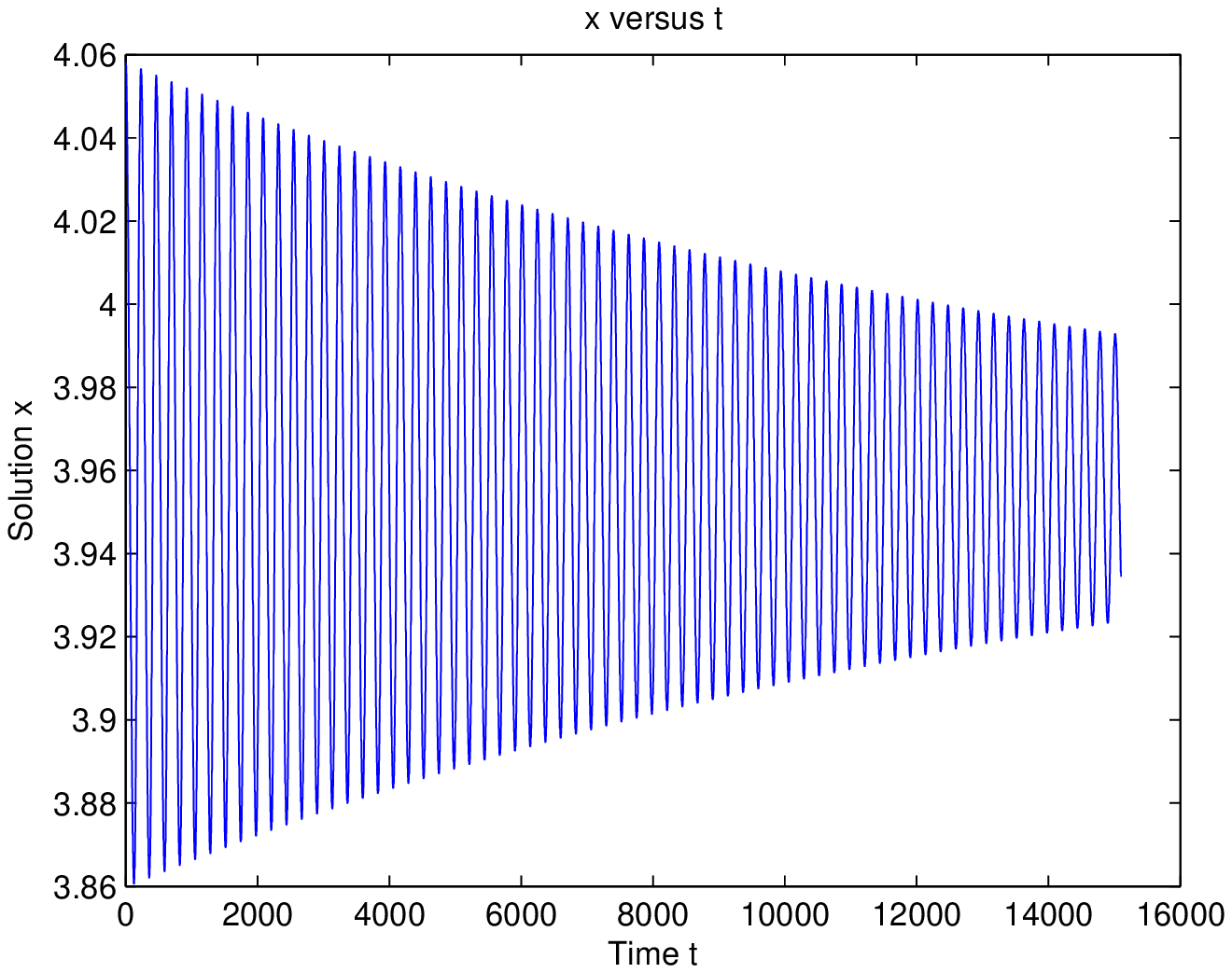}\includegraphics[width=0.49\linewidth]{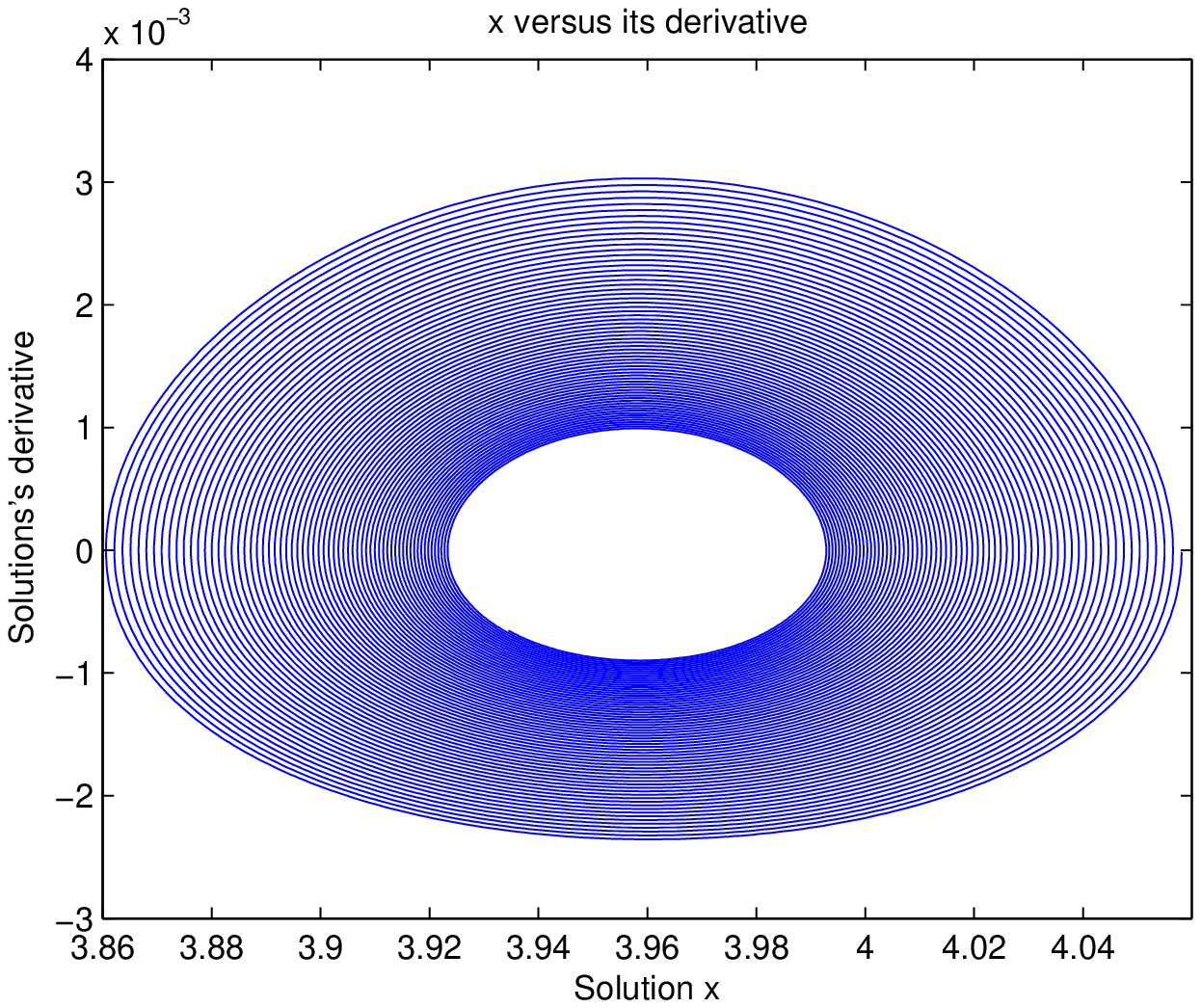}
\caption{\small{Phase portrait for the parameters in the point
$P_3$, for $c=0.1$.}}
\end{figure}

It was very difficult to find a point in the parameter plane,
presenting the behavior corresponding to that of zone 3 of the
bifurcation diagram. We suppose this is so because of the curve $T$
that is deformed and may be very close to the curve $\omega=0$.
However, we found that for the point $P_3 \,(d=0.0015,\,r=7.55),$
the behavior of the solution is the following: for initial functions
close to $x_2$, that is $c\leq 0.42,$ the solution spirals towards
$x_2$, while for $c\geq 0.425$ it is visible that the solution
spirals away from a repulsive limit cycle, having increasing (with
time) amplitude. When time increases, the solution tends to a large
attractive limit cycle (that was previously formed by Hopf
bifurcation when passing from zone 1 to zone 2). This types of
behavior are shown in Figs. 9-13 where we took
$c=0.1,\,c=0.42,\,c=0.425,\,c=0.45 ,\,c=0.6.$ In Figs. 10 and 11 the
repulsive limit cycle is clearly visible, while in Figs. 12 and 13
the exterior attractive cycle is present.
\begin{figure}\centering
\includegraphics[width=0.44\linewidth]{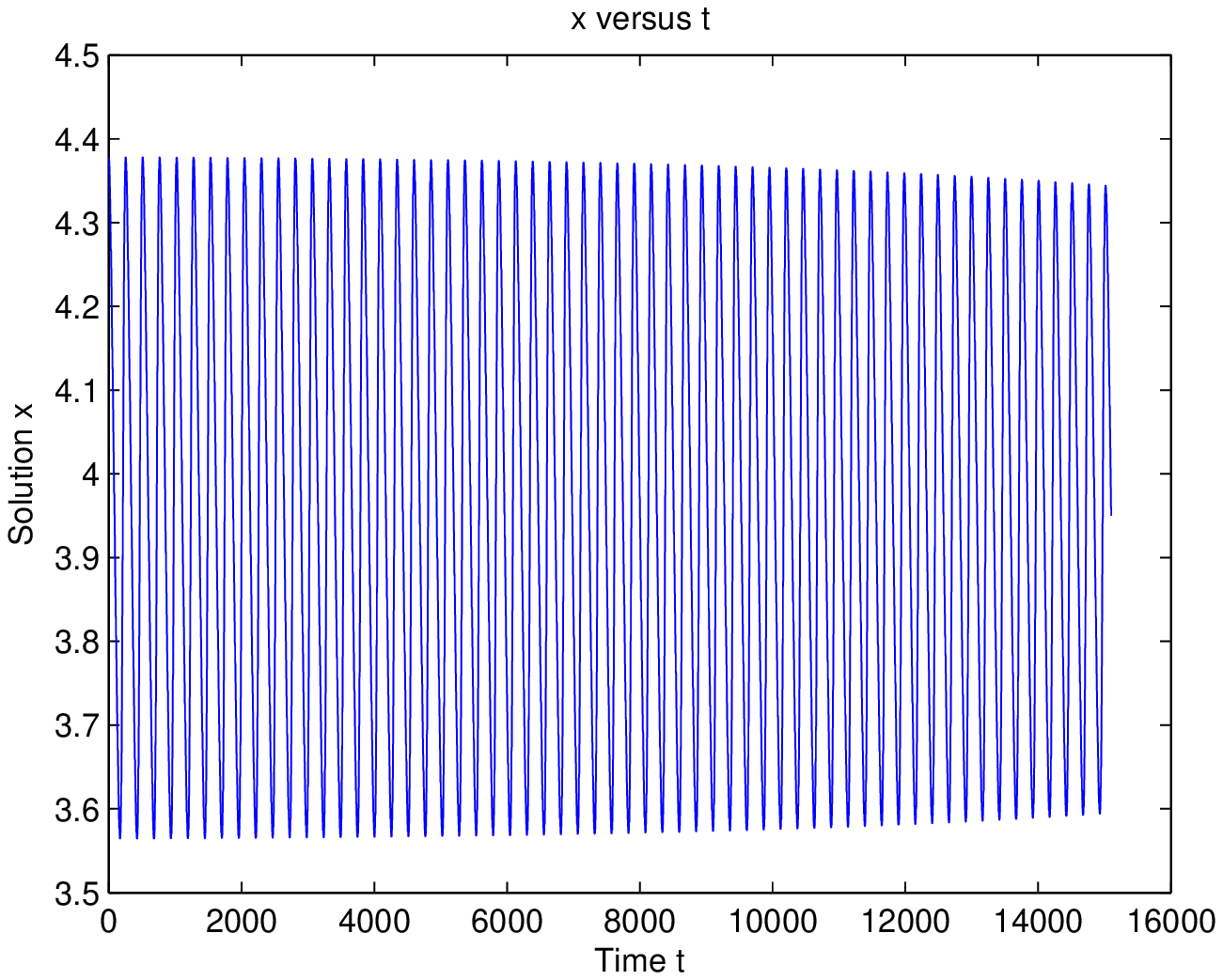}\includegraphics[width=0.49\linewidth]{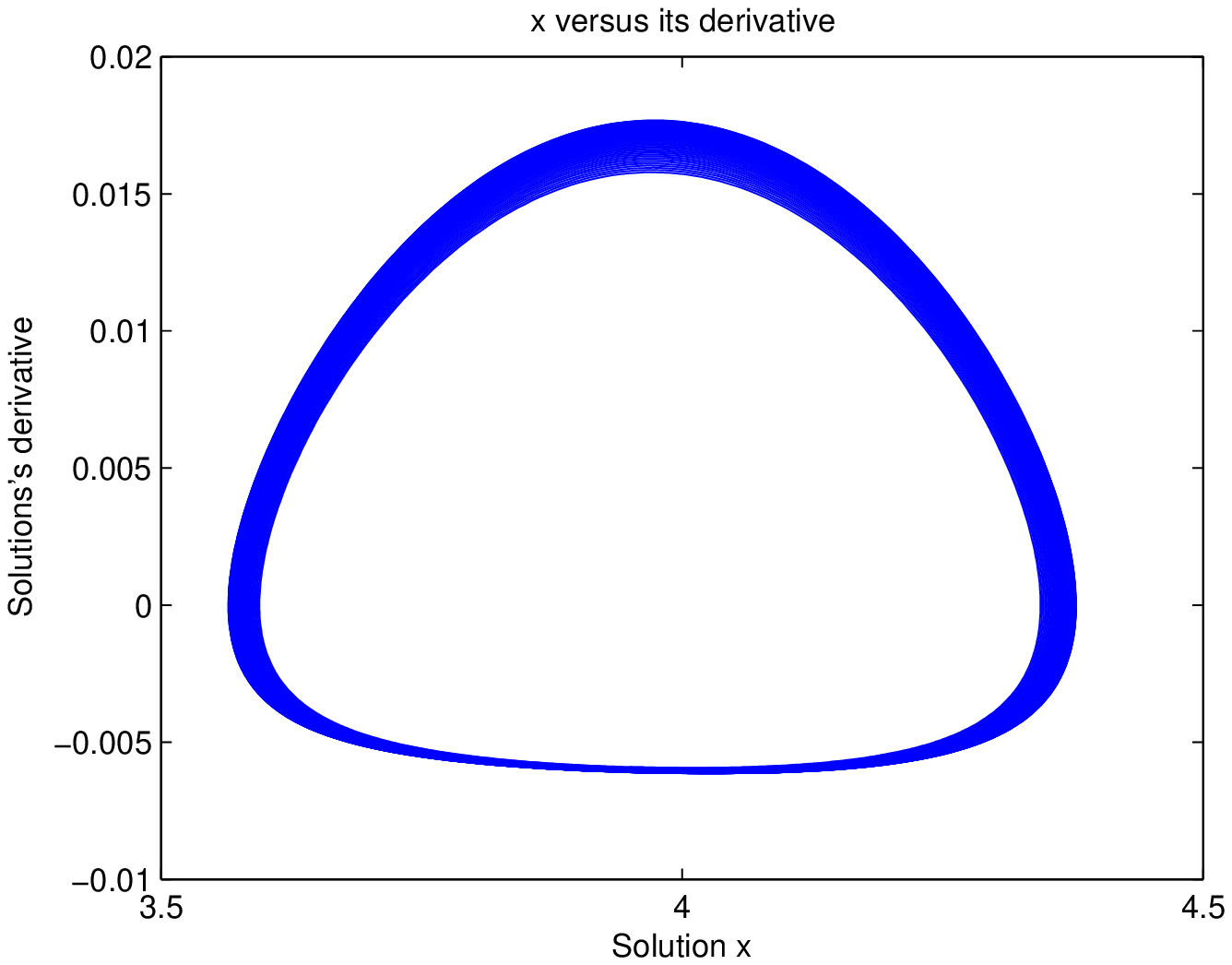}
\caption{\small{Phase portrait for the parameters in the point
$P_3$, for $c=0.42$.}}
\end{figure}

\begin{figure}\centering
\includegraphics[width=0.44\linewidth]{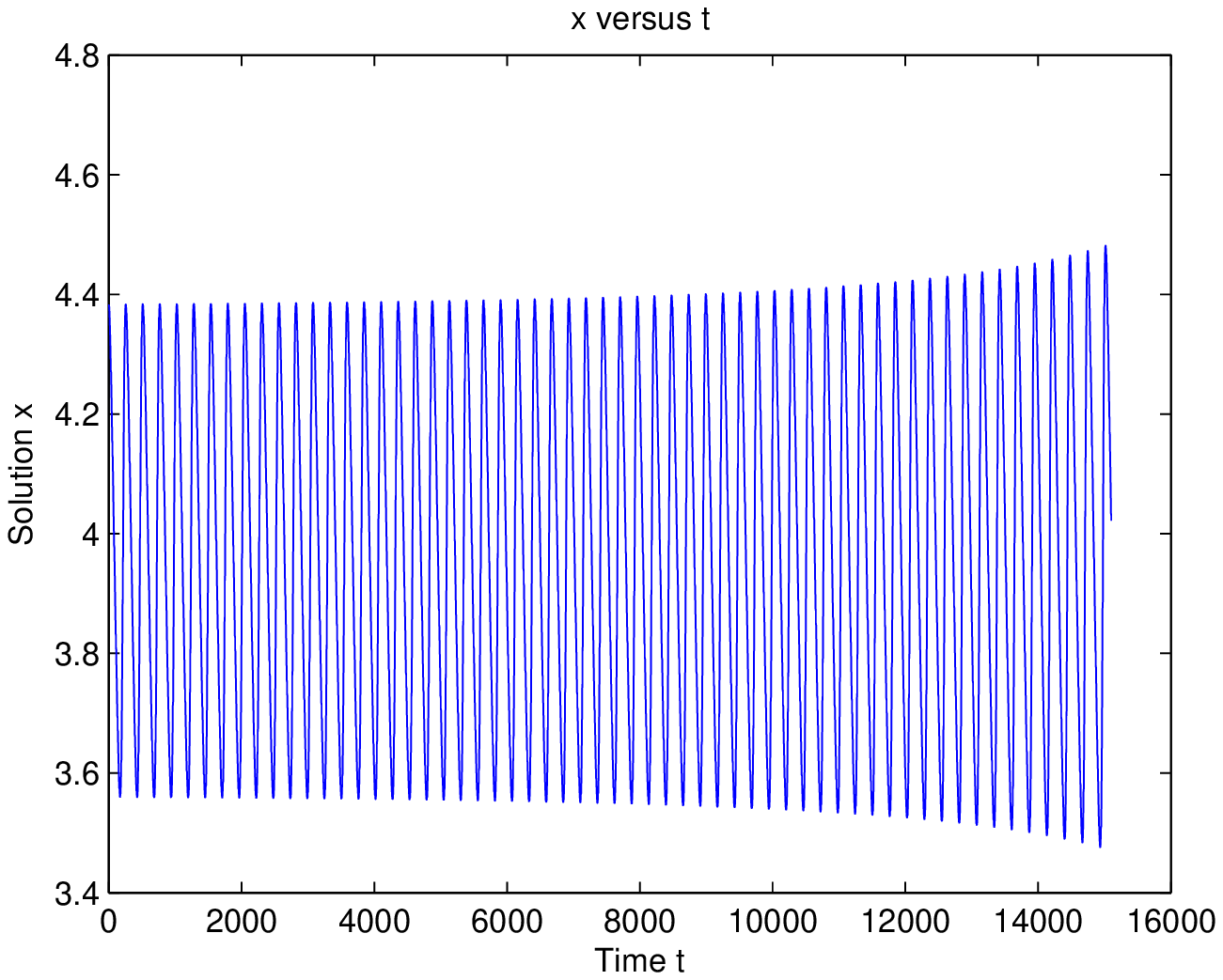}\includegraphics[width=0.49\linewidth]{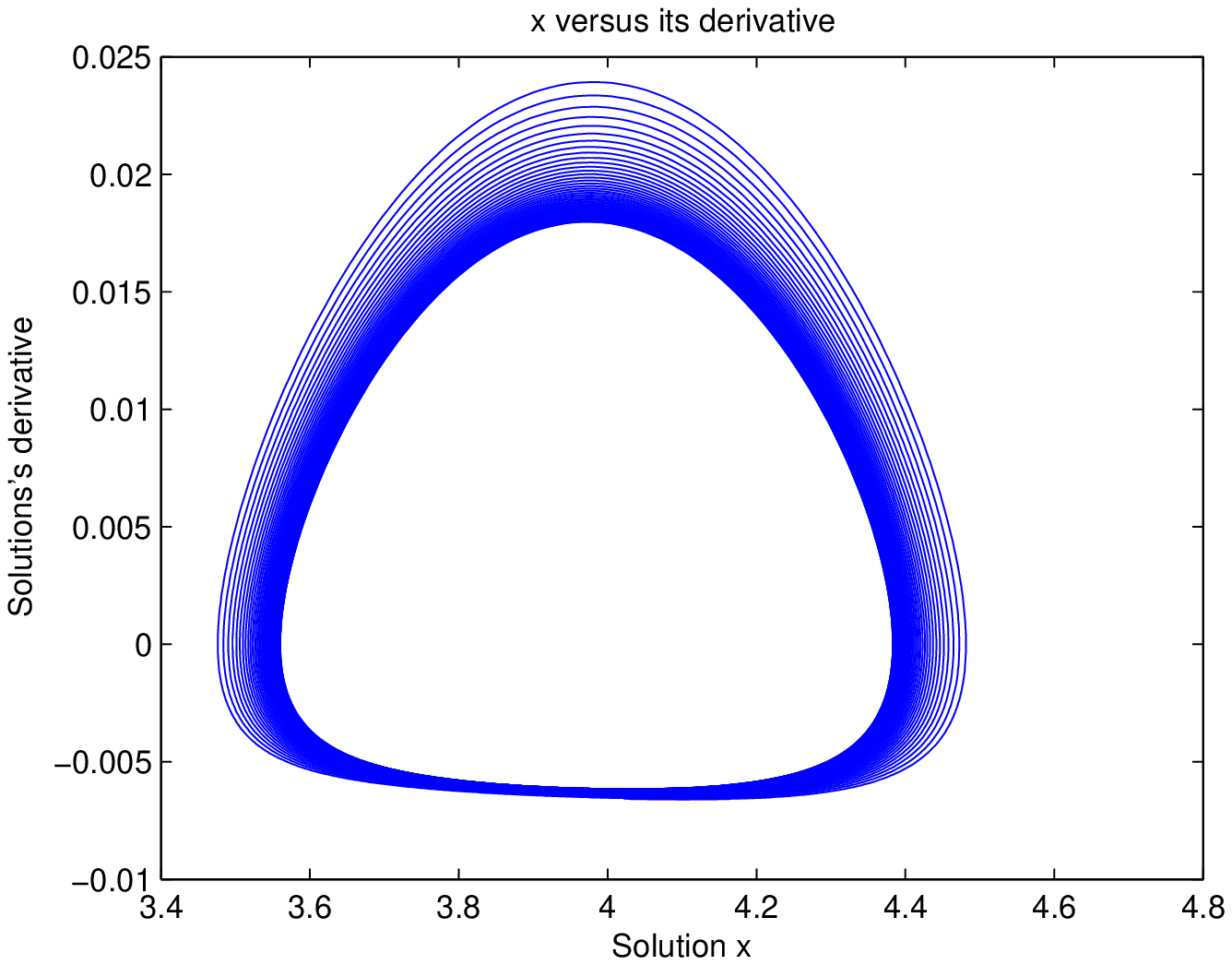}
\caption{\small{Phase portrait for the parameters in the point
$P_3$, for $c=0.425$.}}
\end{figure}

\begin{figure}\centering
\includegraphics[width=0.44\linewidth]{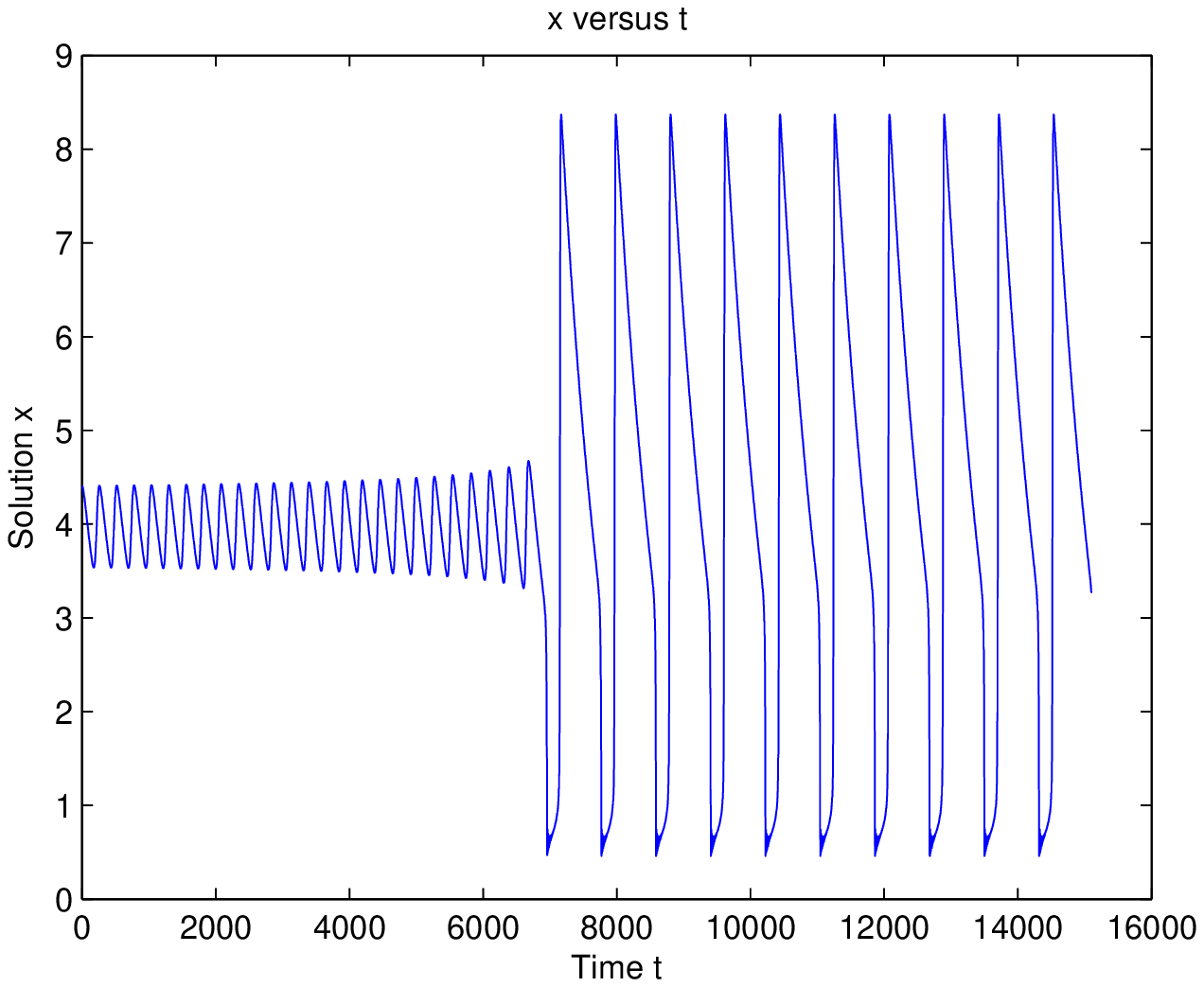}\includegraphics[width=0.49\linewidth]{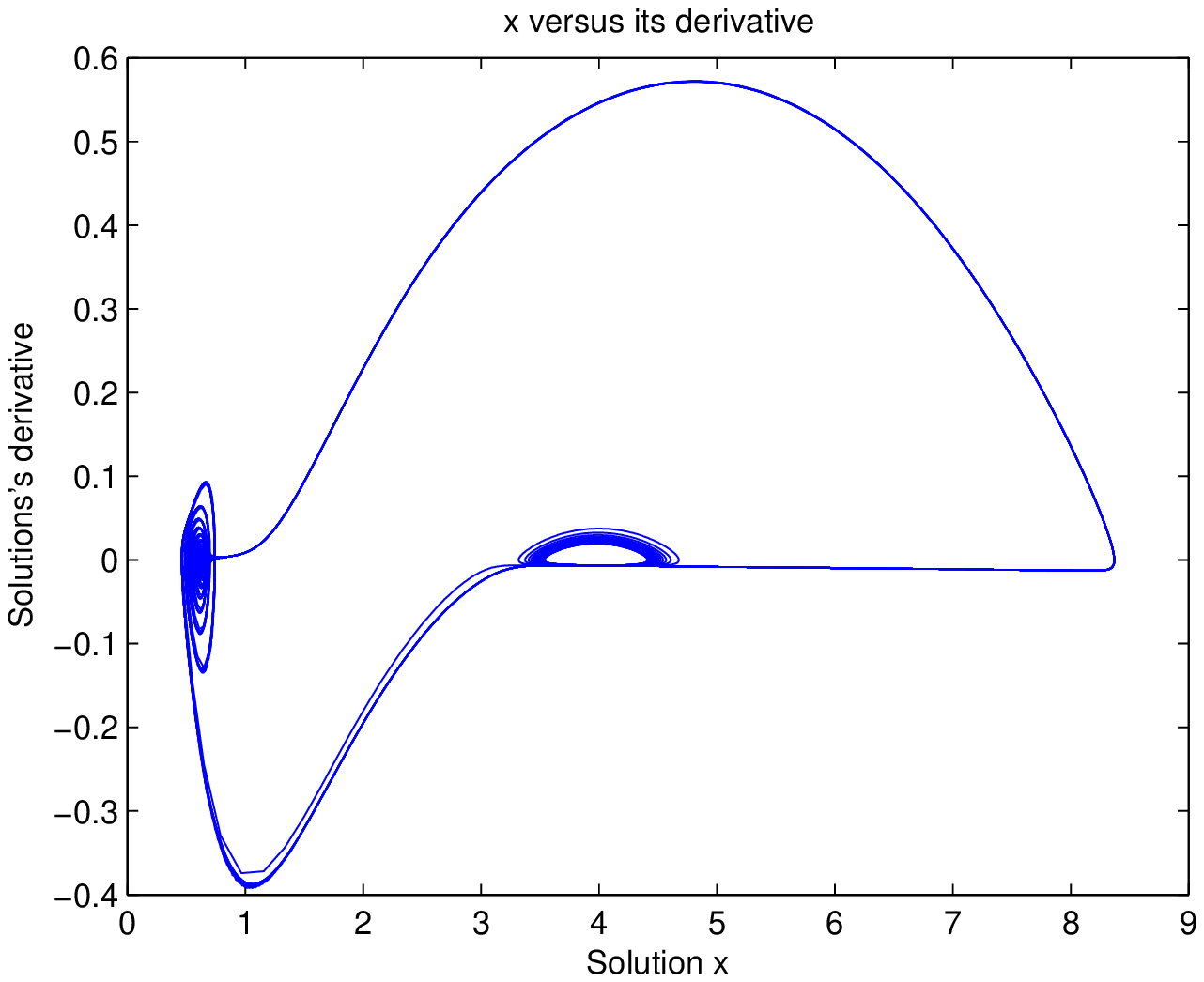}
\caption{\small{Phase portrait for the parameters in the point
$P_3$, for $c=0.45$.}}
\end{figure}
\begin{figure}\centering
\includegraphics[width=0.44\linewidth]{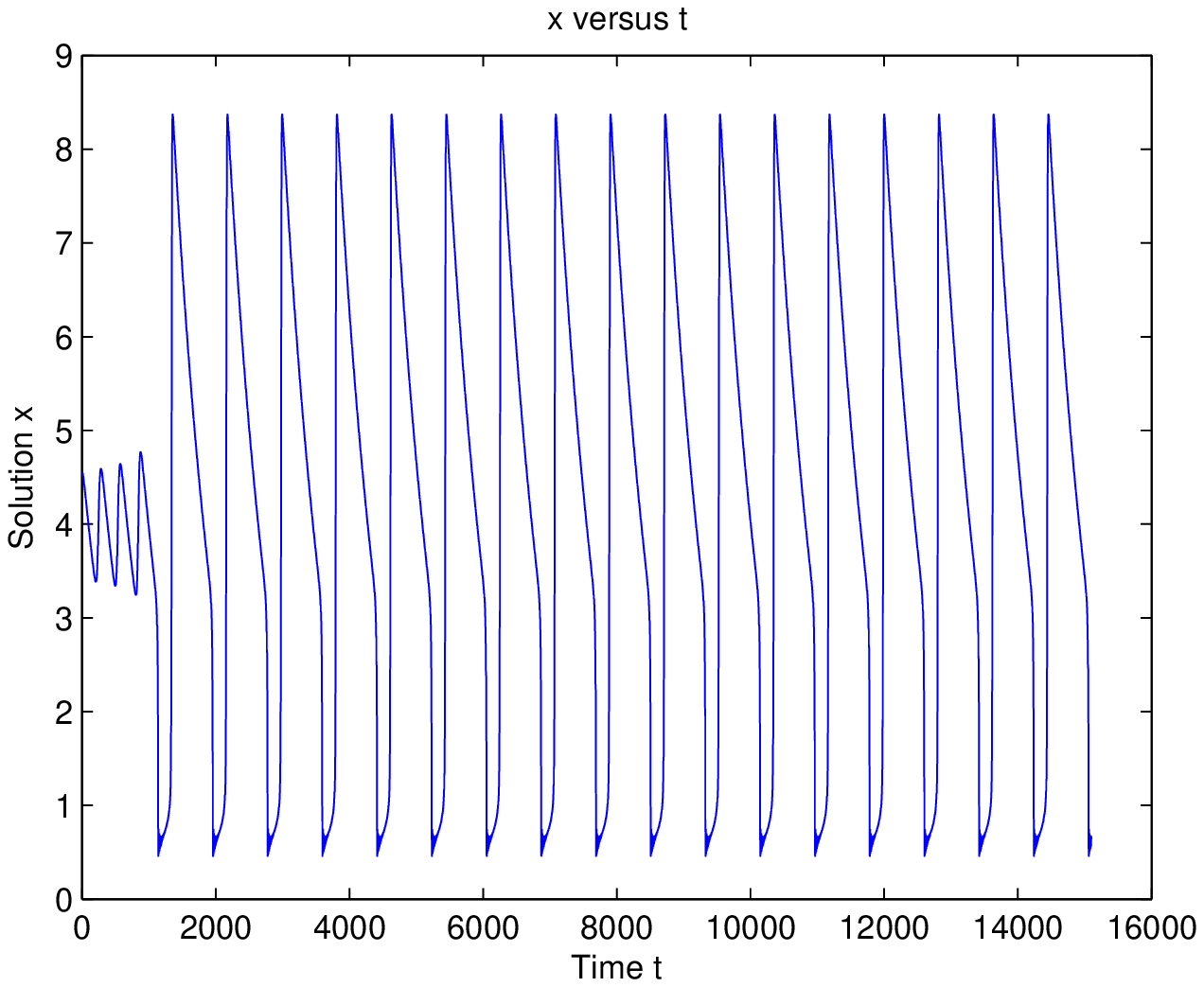}\includegraphics[width=0.49\linewidth]{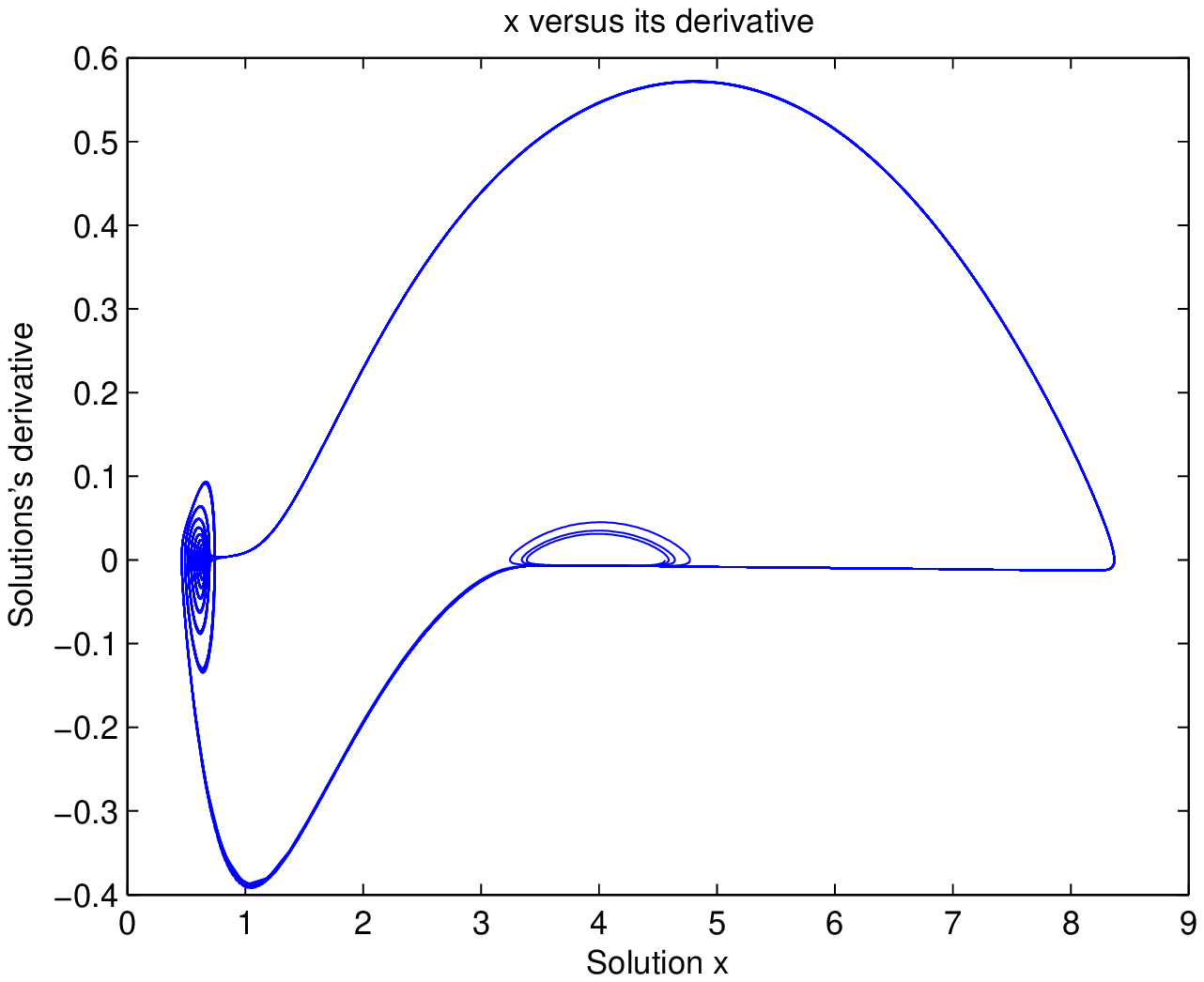}
\caption{\small{Phase portrait for the parameters in the point
$P_3$, for $c=0.6$.}}
\end{figure}

Hence, by numerical integrations we confirmed that the Bautin
bifurcation, predicted by theoretical considerations, actually takes
place, for the considered differential delay equation, in one of the
points with $\l_1=0,\,l_2<0.$

It is important to remark that  a consequence of the Bautin
bifurcation is, for a certain zone in the parameter space,  the
occurrence of a repulsive limit cycle inside an attractive one.
There, the behavior of the solution of eq. (1) strongly depends on
the initial function. If the initial function is close to the
equilibrium point $x_2$, the solution spirals towards $x_2,$ while
if the initial function is ``far'' from $x_2$, it will spiral
towards the exterior stable limit cycle. From the point of view of
the studied illness these two behaviors are quite different and this
shows the importance of being able to control the initial condition.

\end{document}